\documentclass[]{amsart}

\usepackage{amssymb}
\usepackage{amsmath}
\usepackage{amsthm}
\usepackage{amscd}

\usepackage{overpic}

\long\def\comment#1\endcomment{}
\usepackage{color}

\newcommand {\cals}{\mathcal S} 
\newcommand{\base}{\operatorname{base}}

\newcommand {\ft}{\mathfrak t}

\newcommand {\iv}{^{-1}}

\newcommand{\dv}{{\mathrm{div}}}

\newcommand{\Dv}{{\mathrm{Div}}}
\newcommand{\Is}{{\mathrm{Iso}}}

\newcommand{\vol}{{\mathrm{Vol}}}
\newcommand{\fvol}{{\mathrm{FillVol}}}

\newcommand{\Rk}{{\mathrm{Rank}}}

\newcommand {\nn}{\mathcal N} 
\newcommand{\onn}{\overline{{\mathcal N}}}
\newcommand {\calm}{\mathcal M} 

\def\<{\left\langle}
\def\>{\right\rangle}

\def\MM{\mathcal {K}}

\def\CC{{\mathcal C}}

\def\MCG{\mathcal{MCG}}
\def\QQ{{\mathcal Q}}
\newcommand {\pgot}{{\mathfrak p}}
\newcommand {\rgot}{{\mathfrak r}}
\newcommand {\sgot}{{\mathfrak s}}

\newcommand {\cgot}{{\mathfrak c}}

\newcommand {\dgot}{{\mathfrak d}}

\newcommand{\me}{\medskip}

\newcommand {\p}{\mathfrak p} 
\newcommand {\q}{\mathfrak q} 
\newcommand {\pg}{\mathfrak g} 

\newcommand{\fh}{\mathfrak{h}}

\def\MCG{{\mathcal MCG}}

\newtheorem{thm}{Theorem}[section]
\newtheorem{prop}[thm]{Proposition}
\newtheorem{cor}[thm]{Corollary}

\newtheorem{lem}[thm]{Lemma}

\theoremstyle{definition} \newtheorem{defn}[thm]{Definition}

\theoremstyle{remark}
\newtheorem{rmk}[thm]{Remark}

\def\square{\hfill${\vcenter{\vbox{\hrule height.4pt \hbox{\vrule
width.4pt height7pt \kern7pt \vrule width.4pt} \hrule
height.4pt}}}$}

\newcommand{\tsh}[1]{\left\{\kern-.9ex\left\{#1\right\}\kern-.9ex\right\}}

\def\R{{\mathbb R}}

\def\Z{{\mathbb Z}}
\def\N{{\mathbb N}}

\def\dd{{\mathcal D}}

\def\CC{{\mathcal C}}
\def\calc{{\mathcal C}}

\def\V{{\mathcal V}}

\def\P{{\mathcal P}}

\def\nn{{\mathcal N}}

\def\fg{{\mathfrak g}}

\def\fh{{\mathfrak h}}

\def\MCG{\mathcal {MCG}}
\def\MCGS{\mathcal {MCG}(S)}

\def\diam{{\rm{diam}}}

\def\dist{{\rm{dist}}}

\def\<{\langle}
\def\>{\rangle}

\long\def\Restate#1#2#3#4{
\medskip\par\noindent
{\bf #1 \ref{#2} #3} {\it #4}\par\medskip }

\definecolor{darkgreen}{cmyk}{1,0,1,.2}

\newcommand\intersect\cap
\newcommand\infinity\infty
\newcommand\wt\widetilde

\newcommand\inject\hookrightarrow
\newcommand\union\cup

\newcommand{\co}{\colon\thinspace}

\newcommand\join\Lambda
\newcommand\cross\times

\newcommand\lub\vee
\newcommand\glb\wedge

\renewcommand\paragraph[1]{\medskip\textbf{#1} }

\long\def\Restate#1#2{
\medskip\par\noindent
{\bf #1} {\it #2}\par\medskip }

\begin{document}
\title[Divergence for
mapping class groups and CAT(0) groups]{Higher dimensional divergence for
mapping class groups and CAT(0) groups}
\author{Jason Behrstock}
\thanks{The research of the first author was supported
as an Alfred P. Sloan Fellow and by
NSF grant \#DMS-1006219.}
\address{Lehman College and the Graduate Center,
City University of New York,
U.S.A.}
\author{Cornelia Dru\c{t}u}\thanks{The research of the second author was supported in part by
the EPSRC grant ``Geometric and analytic aspects of infinite groups", by the project ANR Blanc ANR-10-BLAN 0116, acronym GGAA, and by the LABEX CEMPI}

\address{Mathematical Institute,
Radcliffe Observatory Quarter,
Oxford OX2 6GG, U. K.}

\subjclass[2010]{{Primary 20F65; Secondary 20F69, 20F38, 20F67}} \keywords{mapping class group, non-positive curvature, isoperimetric function, divergence, rank}

\date{\today}

\begin{abstract} In this paper we investigate the higher dimensional
divergence functions of mapping class groups of surfaces and of
$CAT(0)$--groups.  We show that, for mapping class groups of surfaces,
these functions exhibit phase transitions at the rank (as measured by
$3\cdot\mbox{genus}+\mbox{number of punctures}-3$).  We also provide
inductive constructions of $CAT(0)$--spaces with co-compact group
actions, for which the divergence below the rank is (exactly) a
polynomial function of our choice, with degree arbitrarily large
compared to the dimension.
\end{abstract}

\maketitle

\section{Introduction} 

In this paper we study \emph{higher dimensional divergence functions},
a type of filling function which is particularly significant in spaces with an
interesting geometry on the boundary at infinity, such as spaces of
non-positive curvature.  In some sense, these functions describe the
spread of geodesics, and the filling close to the boundary at
infinity.  They provide a powerful quasi-isometric invariant with 
which to study the large-scale
geometry of a space or a group.

Given a fixed point $x_0$ in a metric space $X$, the $k$--\emph{dimensional divergence function} (or the $k$--\emph{divergence}) \emph{in} $r$, roughly speaking, measures the minimal filling of $k$--dimensional spheres avoiding the ball
$B(x_0,r)$, and of area at most $Ar^k$, by $(k+1)$--dimensional balls that are outside the open ball $B(x_0, \lambda r)$. The choices of the fixed
center,
$x_0$,  and of the parameters $A>0$ and $\lambda \in (0,1)$ do not
affect the order of the $k$--divergence. 
Versions of divergence functions were defined for non-positively
curved manifolds in
\cite{BradyFarb:div} (there, spheres, and the balls filling
them, are Lipschitz maps of the corresponding Euclidean spheres and
balls into the given manifold, and volumes are computed using the fact
that such maps are differentiable almost everywhere) and for
$CAT(0)$--spaces and integral currents in
\cite{Wenger:div}. 

In the present paper,
we work in the setting of finitely generated groups, so we will
use classifying spaces and a version of divergence defined using
combinatorial filling functions. The $k$--divergence functions, for 
$k$ at least $1$, are also called \emph{higher dimensional} divergence
functions, by contrast with the $0$--divergence functions which by
their very nature are less of isoperimetric filling functions, and
more of distortion functions of exteriors of large balls.

We investigate the behavior of higher dimensional divergence functions
for $CAT(0)$--groups (i.e., groups acting properly discontinuously
cocompactly on a $CAT(0)$--space) and for mapping class groups
$\MCG(S)$ of surfaces $S$ (i.e., groups of isotopy classes of
homeomorphisms of surfaces $S$).  It is known that mapping class
groups have a large scale geometry that is $CAT(0)$, in the sense that 
all their asymptotic cones are bi-Lipschitz equivalent to
$CAT(0)$--spaces with a uniform bi-Lipschitz constant
\cite{BehrstockDrutuSapir:mcg,Bowditch:cat0median}; other 
similiarities have been exposed in \cite{BehrstockHagenSisto:HHS_I}.
On the other hand, they are not actual $CAT(0)$--groups as they cannot
act geometrically on any $CAT(0)$--space when $S$ has genus at least
three, or genus $2$ and at least one puncture
\cite{KapovichLeeb:actions}; when $S$ is a closed surface of genus
$g$, $\MCG(S)$ does not admit any action by semisimple isometries on
complete $CAT(0)$--metric spaces of dimension less than $g$
\cite{Bridson:mcgCAT}.
The present paper provides further similiarities between mapping
class groups and $CAT(0)$--groups, in particular we introduce some methods
to provide lower bounds for higher dimensional divergence that are
common to the two types of groups.

A central notion in non-positively curved geometry is that of
\emph{rank} which is often taken as the largest
dimension at which Euclidean behavior can still be detected, the 
meaning of ``behavior'' can vary as we see below; in a 
sense the rank gives a dimension above which
hyperbolic features are bound to appear.  There
are several ways to define the rank, and the connections between the
different
definitions are not fully understood (see for instance
\cite[$\S 6.B_2$]{Gromov:Asymptotic}, where seven possible definitions of rank are listed). To begin with, there exists the notion of \emph{(quasi-)flat rank}, i.e., the maximal dimension of a Euclidean
space that can be (quasi-)isometrically embedded in the given space. In the case of a locally compact complete simply connected $CAT(0)$--space that admits a cocompact action by isometries, the quasi-flat rank equals the flat rank
\cite{AndersonSchroeder,Gromov:Asymptotic,Kleiner:lengthsp}; thus this notion of rank has the two-fold advantage of being consistent with the classical notion of rank in symmetric space, and of being easily transferred to groups acting properly cocompactly on $CAT(0)$--spaces, \textit{via} its quasi-flat version.
The quasi-flat rank of a mapping class group $\MCG (S)$ equals the maximal rank of
a free abelian subgroup of
$\MCG (S)$, i.e., the \emph{complexity} of $S$, as
measured by
$\xi (S) = 3g +p -3$, where
$g$ is the genus and $p$ is the number of boundary components of $S$.
(We do not require that mapping classes fix the ``boundary'', so
formally each boundary component should be considered as a ``puncture'').
That the algebraic rank is equal to $\xi(S)$ was proven by
\cite{BirmanLubotzkyMcCarthy}, the equality between the algebraic rank and
the quasi-flat rank was proven in \cite{BehrstockMinsky:rankconj,
Hamenstadt:qirigidity}.

Another notion of rank is defined using isoperimetric
functions. One considers the rank in this sense to be
the maximal dimension for which these functions are of
the same order as those in Euclidean spaces.

Finally, a third notion of rank is obtained by considering the
dimension in which there is a phase
transition for the higher dimensional divergence functions.

For mapping class groups, we prove that these ranks all agree. 
That the isoperimetric rank equals the quasi-flat rank follows from Theorems
\ref{thm:genus3} and \ref{thm:genus012}.
A main result of this paper is on the higher dimensional
divergence functions of the mapping
class group; one aspect of this theorem is that of displaying a ``phase
transition'' for these
functions, when the dimension equals the quasi-flat rank, whence 
providing equality of the first two notions of rank with the third.

\Restate{Theorems \ref{thm:divOverRk} and \ref{thm:divsubrank}}{
Let $S$ be a compact orientable surface
and let
$\Dv_k$ be the $k$--dimensional divergence of $\MCG (S)$.

If $k< \xi (S)$ then $\Dv_{k}\succeq x^{k+2}$.

If $k\geq \xi (S)$ then $\Dv_k (x) = o(x^{k+1})$; further, if $S$
is either genus $0$ or $1$, or genus $2$ with empty boundary, then
$\Dv_k (x) \asymp x^{k}$.
}

The relation of one function  being 
``asymptotically bounded'' by another, $\succeq$, or two functions 
being asymptotically equivalent, $\asymp$, are made precise in 
Section \ref{sect:prelim}.

A special case of the above result is that it establishes that the $1$-dimensional divergence is at least cubic for surfaces with $\xi(S) \geq 2$, while it is proved in \cite{ABDDY:Pushing} that for such
surfaces the
$1$--dimensional divergence is at most quartic. It would be interesting to know the exact order of the $1$-dimensional divergence. Besides the estimate in \cite{ABDDY:Pushing}, the only other 
previously known result about divergence in mapping class groups is that the $0$--divergence is quadratic \cite{Behrstock:asymptotic}.

In the case of $CAT(0)$--spaces, the known facts on divergence are that: for a symmetric space $X$ of
non-compact type, $\Dv_k $ grows exponentially when $k=\Rk (X)-1$, see
\cite{BradyFarb:div,Leuzinger:divergence}, while when $k\geq \Rk (X)$
the divergence satisfies $\Dv_k = O( x^{k+1})$, see \cite{Hindawi:div}.
For a
cocompact $CAT(0)$--space $X$ and a homological
version of the divergence defined in terms of integral currents, it has been established \cite{Wenger:div} that if
$k=\Rk (X)-1$ then
$\Dv_k \succeq x^{k+2}$, while if $k\geq \Rk (X)$ then $\Dv_k \preceq
x^{k+1}$.

The  picture becomes unexpectedly complicated when it comes to divergence below the rank.  The following theorem yields surprising examples in the context of $CAT(0)$--spaces.

\Restate{Theorem \ref{thm:divmcg}}{
For every positive integers $r$ and $n$ there exist universal covers
of compact CAT(0)--spaces with flat rank $2r$, and such that the
$(r-1)$--dimensional divergence satisfies $\Dv_{r-1}\asymp 
x^{r+n}$.
}

\medskip

In our study of divergence, we develop a method to obtain lower bounds for higher divergence from lower bounds on zero--divergence. We use this method both for mapping class groups and for $CAT(0)$--spaces. 

For the inductive construction of $CAT(0)$--spaces with pathological
divergence behaviour described in the above theorem, the sharp
estimate from above of the divergence is obtained by using results
from \cite{BD:HigherFilling1} (see Theorem \ref{thm:DivRound}), which 
allow us to restrict our study to round spheres (i.e., $k$--dimensional
spheres with diameter controlled by the $k$-th root of the volume),
and a careful study of the structure of round spheres in the
$CAT(0)$--spaces that we construct.

\medskip

The plan of the paper is as follows.  In Section \ref{sect:prelim0} we
recall some basic notions and establish notation which we will use in
the paper, in particular some combinatorial
terminology and results on mapping class groups. Also, we recall 
background from \cite{BD:HigherFilling1} on combinatorial
formulations of isoperimetric and divergence functions and 
their behavior in the presence of a combing.  Section \ref{sec:divmcg} contains the theorem describing the behavior of divergence functions in mapping class groups. Section \ref{sec:divcat}
contains the proof of Theorem \ref{thm:divmcg}.

\subsection*{Acknowledgements}
The authors thank Bruce Kleiner, Enrico Leuzinger and Stefan Wenger for useful comments on an earlier version of this paper. Part of the work on this paper was carried during visits of the first author to the Mathematical Institute in Oxford. We thank the institute for its hospitality.

\section{Preliminaries}\label{sect:prelim0}

\subsection{General terminology}\label{sect:prelim}

Given two functions $f,g$ which both map $\R_+$ to itself, we write $f\preceq_{C,k} g$ for some constant $C\geq 1$ and integer $k\geq 1$, if
$$f(x) \le Cg(Cx+C)+Cx^{k}+C\; \; \mbox{ for all }x\in \R_+.
$$
 
We write $f\asymp_{C,k} g$ if and only if
$f\preceq_{C,k} g$ and $g\preceq_{C,k} f$. 
Two functions $\R_+\to \R_+$ are
said to be $k$--\emph{asymptotically equal} if there exists $C\geq 1$ s.t.
$f\asymp_{C,k} g$.  This is an equivalence relation.  

When at least one of the two functions $f,g$ involved in the relations above is an $n$--dimensional isoperimetric or divergence function, we automatically consider only relations where $k=n$, therefore $k$ will no longer appear in the subscript of the relation. When irrelevant, we
do not mention the constant $C$ either and likewise remove the corresponding subscript.

We use the notations $f=O(g)$ and $f=o(g)$, for $f$ and $g$ real-valued functions of one real variable, with their standard meaning.

\me

In a metric space $(X, \dist )$, the {\em open $R$--neighborhood} of
a subset $A$, i.e. $\{x\in X:
\dist (x, A)<R\}$, is denoted by $\nn_R(A)$.  In particular, if $A=\{a\}$ then $\nn_R(A)=B(a,R)$
is the open $R$--ball centered at $a$. We use the notation $\onn_R (A)$ and $\bar{B}(a,R)$ to designate the
corresponding {\em closed neighborhood} and {\em closed ball}
defined by non-strict inequalities. We make the convention that $B(a,R)$ and $\bar{B}(a,R)$ are the empty
set for $R<0$ and any $a\in X$.

\me

Fix two constants $L\geq 1$ and $C\geq 0$. A map $\q\co Y\to X$ is
said to be
\begin{itemize}
  \item  $(L,C)$--\emph{quasi-Lipschitz} if
$$
\dist (\q (y),\q (y'))\leq L\dist (y,y')+C,
\hbox{ for all } y,y' \in Y\, ;
$$
  \item  an $(L,C)$--\emph{quasi-isometric embedding} if moreover
  $$
\dist (\q (y),\q (y')) \geq \frac{1}{L}(y,y')-C\hbox{ for all } y,y' \in Y\, ;
$$
  \item an $(L,C)$-\emph{quasi-isometry} if it is an $(L,C)$--quasi-isometric
embedding $\q\co Y\to X$ satisfying the additional assumption that $X\subset \nn_C(\q
(Y))$.
  \item an $(L,C)$--{\em quasi-geodesic} if it is an $(L,C)$--quasi-isometric
embedding defined on an interval of the real
line;

\item a \emph{bi-infinite $(L,C)$--quasi-geodesic} when defined on the entire real line.
\end{itemize}

In the last two cases the terminology is extended to the image of $\q$. When the constants $L,C$ are
irrelevant they are not mentioned.

We call $(L,0)$--quasi-isometries (quasi-geodesics)
$L$--\emph{bi-Lipschitz maps (paths)}.
If an $(L,C)$--quasi-geodesic $\q$ is $L$--Lipschitz then $\q$ is called an {\em $(L,C)$--almost geodesic}.
Every
$(L,C)$-quasi-geodesic in a geodesic metric space is at bounded (in
terms of $L, C$) distance from an $(L+C,C)$--almost geodesic with the
same end points, see e.g. \cite[Proposition 8.3.4]{Buragos:course}.
Therefore, without loss of generality, we assume in
this text that all quasi-geodesics are in fact almost geodesics, in
particular that they are continuous.

\me

Given two subsets $A,B\subset \R$,
a map $f\co A\to B$ is said to be
\emph{coarsely increasing} if there exists a constant $D$ such that
for each $a,b$ in $A$ satisfying $a+D<b$, we have that $f(a)\leq f(b)$. Similarly, we define
\emph{coarsely decreasing} and \emph{coarsely monotonic} maps. A
map between quasi-geodesics is coarsely monotonic if it defines a
coarsely monotonic map between suitable nets in their domain.

\me

A metric space is called
\begin{itemize}
  \item \emph{proper} if all its closed balls are compact;
  \item \emph{cocompact} if there exists a compact subset $K$ in $X$ such that all the translations of $K$ by isometries of $X$ cover $X$;
  \item \emph{periodic} if it is geodesic
and for fixed constants $L\geq 1$ and $C\geq 0$ the image of
some fixed ball under $(L,C)$--quasi-isometries of $X$ covers $X$;
  \item a \emph{Hadamard space} if $X$ is geodesic, complete, simply connected and satisfies the CAT(0)
condition;
  \item  a \emph{Hadamard manifold} if moreover $X$ is a
smooth Riemannian manifold.
\end{itemize}

\subsection{Combinatorial terminology}\label{subsec:combi}

We use the standard terminology related to simplicial complexes as it appears in \cite{Hatcher:book}. When we speak of
simplicial complexes in what follows, we mean their topological realisation most of the time. Throughout the paper, we assume that all simplicial complexes are connected. We consider every simplicial complex endowed with a ``large scale metric
structure'' defined by assuming that all edges have length one and taking the
shortest path metric on the $1$-skeleton. 

We say that a simplicial complex $X$ has a \emph{bounded
$(L,C)$--quasi-geodesic combing}, where $L\geq 1$ and $C\geq 0$, if
for every $x\in X^{(1)}$ there exists a way to assign to every element
$y\in X^{(1)}$ an $(L,C)$--quasi-geodesic $\q_{xy}$ connecting $y$ to
$x$ in $X^{(1)}$, such that
$$
\dist (\q_{xy}(i), \q_{xa}(i))\leq L\dist (y,a)+L\, ,
$$ for all $x, y,a\in X^{(1)}$ and $i\in \R$.  Here the quasi-geodesics are assumed
to be extended to $\R$ by constant maps.

 Given an $n$--dimensional simplicial complex $\calc$, we call
  the closed simplices of dimension $n$ the
  \emph{chambers} of $\calc$. 

Given a simplicial map $f\co X\to Y$, where $X,Y$ are
simplicial complexes, $X$ of dimension $n$, we call $f$--\emph{non-collapsed chambers in $X$} the chambers whose images by $f$ stay of dimension $n$. We denote by $X_\vol$ the set of $f$--non-collapsed chambers. We define the \emph{volume of }$f$ to be the (possibly infinite) cardinality of $X_\vol$.

\me

Recall that a group $G$ is \emph{of type} $\mathcal{F}_k$ if it admits
an Eilenberg-MacLane space $K(G,1)$ whose $k$-skeleton is finite.

\begin{prop}[\cite{AWP}, Proposition 2]\label{prop:typefn}
If a group $G$ acts cellularly on a CW-complex $X$, with finite
stabilizers of points and such that $X^{(1)}/G$ is finite then $G$ is
finitely generated and quasi-isometric to $X$.  Moreover, if $X$ is
$n$-connected and $X^{(n+1)}/G$ is finite then $G$ is of type
$\mathcal{F}_{n+1}$.
\end{prop}

Conversely, it is easily seen that for a group of type
$\mathcal{F}_{n+1}$ one can define an $(n+1)$-dimensional
$n$-connected simplicial complex $X$ on which $G$ acts
properly discontinuously by simplicial isomorphisms, with trivial
stabilizers of vertices, such that $X/G$ has finitely many cells.  Any
two such complexes $X,Y$ are quasi-isometric, and the quasi-isometry,
which can initially be seen as a bi-Lipschitz map between two subsets
of vertices, can be easily extended to a simplicial map $X\to Y\, $
\cite[Lemma 12]{AWP}.

A group is \emph{of type} $\mathcal{F}_\infty$ if and only if it is of
type $\mathcal{F}_k$ for every $k\in \N\, $.
It was proven in \cite[Theorem 10.2.6]{ECHLPT} that every combable group is
of type $\mathcal{F}_\infty$.

\subsection{Higher dimensional filling functions}\label{sect:higherIso}

Since a finite
$(n+1)$--presentation of a group composed only of simplices can always
be found, it suffices to restrict to simplicial complexes when discussing filling problems. Throughout the paper, we make use of the terminology related to the $n$--dimen\-sio\-nal filling
functions in the simplicial setting as described in full detail in \cite{BD:HigherFilling1}. We briefly recall a number of relevant notions and results.

For the rest of the paper, we assume that all the simplicial complexes that we consider are universal covers of compact simplicial complexes. We also assume that, in the study of filling functions up to dimension $n$, the simplicial complexes considered have dimension ar least $n+1$ and are
$n$--connected.

When we speak of manifolds we always mean
\emph{manifolds with a simplicial-complex structure}.

We denote by $V$ an arbitrary $m$--dimensional
connected compact sub-manifold of $\R^{m}\, $, where $m\geq 2$ is an
integer and $V$ is smooth or piecewise linear, and with boundary.  We
denote its boundary by $\partial V$. Unless otherwise stated, the standing assumption is that $\partial V$ is connected.

Given $V$ as above, a \emph{domain modelled on $V$ in $X$} (also called a $V$--\emph{domain}) is a simplicial map $\frak d$ of $\dd$ to $X^{(m)}$, where $\dd$ is a
simplicial structure on $V$.  When the manifold $V$ is irrelevant we
simply call $\frak d$ a \emph{domain of dimension $m$} (somewhat
abusively, since it might have its entire image inside $X^{(m-1)}$);
we also abuse notation by using $\frak d$ to denote both the map and its image.

A \emph{hypersurface in $X$ modelled on $\partial V$} (also called a $\partial V$--\emph{hypersurface}) is a simplicial map  $\fh$ of
$\calm$ to $X^{(m-1)}$, where $\calm$ is a simplicial structure of the boundary
  $\partial V$.  Again, we abuse notation by letting $\fh$ also denote the image of the above
map, and we also call both $\fh$ and its image a \emph{hypersurface of dimension $m-1$}.

According to the terminology introduced previously, $\dd_\vol$, respectively $\calm_\vol$, is the set of $\frak d$--non-collapsed chambers (respectively $\fh$--non-collapsed chambers). The \emph{volume} of $\frak d$ (respectively $\fh$) is the cardinality of $\dd_\vol$, respectively $\calm_\vol$. 

We sometimes say that the domain $\frak d$ is a $V$--{\emph{domain}}, and $\fh$ is a $\partial V$--{\emph{hypersurface}}. When $V$ is a closed ball in $\R^m$, we call $\frak d$ an $m$--\emph{dimensional ball} and $\fh$ an $(m-1)$--\emph{dimensional sphere}.

\begin{defn}\label{def:round}
A $k$--dimensional hypersurface
 $\fh$ is called
$\eta$--\emph{round} for a constant $\eta >0$ if $\diam (\fh) \leq \eta
\vol (\fh )^{\frac{1}{k}}\, $.
\end{defn}

We say that a domain $\frak d$ \emph{fills} a hypersurface $\fh$ if
this pair corresponds to a $(k+1)$--dimensional connected compact
smooth sub-manifold with boundary $V$ in $\R^{k+1}\, $ satisfying $\dd \cap
\partial V = \calm$ and $\frak d|_\calm = \fh$,
possibly after pre-composing $\fh$ with a simplicial equivalence of
$\calm$.

The \emph{filling volume} of the hypersurface $\fh$,
$\fvol (\fh )$, is the minimum of all the volumes of
domains filling $\fh$.  If no domain filling $\fh$ exists then we set
$\fvol (\fh ) = \infty$.

\begin{defn}\label{def:isopFct}
The $k$\emph{--th isoperimetric function}, also known as the  $k$\emph{--th 
filling function}, of a simplicial complex $X$ is the
function $\Is_k \colon \R_+^* \to \R_+ \cup \{ \infty \}$ such that
$\Is_k (x)$ is the supremum of the filling volume $\fvol (\fh )$ over
all $k$--dimensional spheres $\fh$ of volume at most $Ax^{k}$.
\end{defn}

In what follows the constant $A>0$ from Definition
  \ref{def:isopFct} is fixed, but not made explicit.  Two filling functions corresponding to two different
  values of $A$ are equivalent in the sense of the relation 
  $\asymp$.

The isoperimetric function has a more general version, using instead of the sphere and
its filling with a ball, a hypersurface and its filling with a domain,
both modelled on a $(k+1)$--dimensional submanifold with boundary $V$
in $\R^{k+1}$.  We then define as above the filling function, denoted $\Is_V$.

According to \cite[Theorem 1, Corollary 3]{AWP}, if two simplicial complexes, $X_1$ and $X_2$, of dimension at least $k+1$, are quasi-isometric, then $ \Is_V^{X_1} \asymp \Is_V^{X_2}\, $ for every domain $V$. Therefore the corresponding isoperimetric functions $\Is_{k}$ and $ \Is_V$ are well defined, up to the equivalence relation $\asymp$, for every group  of type $\mathcal{F}_{n+1}$.

\medskip

In \cite{BD:HigherFilling1} we proved that under certain conditions, which are satisfied in the presence of a bounded quasi-geodesic combing, an arbitrary sphere has a
partition into round spheres, such that the sum of
the volumes of the spheres in the partition is bounded by a multiple
of the volume of the initial sphere.  Here a \emph{partition} 
consists of the following: finitely many spheres (more generally,
hypersurfaces) $\fh_1,\dots, \fh_n$ compose a partition of a sphere
$\fh$ if by filling all of them one obtains a ball filling $\fh$ (see
\cite{BD:HigherFilling1}).  The hypersurfaces $\fh_1,\dots,
\fh_n$ are called \emph{contours of the partition}.

\begin{thm}[\cite{BD:HigherFilling1}]\label{thm:rounddec0}
Let $X$ be a simplicial complex with a bounded quasi-geodesic combing, and let $k\geq 2$ be an integer.  

Then for every $\varepsilon >0$ there exists a constant $\eta >0$ such that every $k$--dimensional sphere $\fh$ has a
partition with contours $\fh_1,...\fh_n$ that are $\eta$--round spheres, and contours $\rgot_1,\dots ,\rgot_m$ that are hypersurfaces of volume and filling volume zero such that
\begin{enumerate}

\item\label{rdec1} $\sum_{i=1}^n \vol (\fh_i ) \leq 2\cdot 6^{k+1} \vol (\fh )\, .$

\medskip

\item\label{rdec2} $\fh_1,...\fh_n$ and $\rgot_1,\dots ,\rgot_m$ are  contained in the tubular neighborhood $\nn_R (\fh )$, where $R= \varepsilon \vol (\fh )^{1/k}$.
\end{enumerate}

\end{thm}

A consequence of the above is that isoperimetric functions are bounded by the corresponding Euclidean isoperimetric functions. 

\begin{cor}[The Federer-Fleming inequality for
groups; \cite{BD:HigherFilling1}]\label{cor:rounddec}
Assume that the simplicial complex $X$ has a bounded quasi-geodesic combing.  Then for every $k\geq 1$, $\Is_k (x)
\preceq x^{k+1}$.  Moreover for $k=2$ the supremum of
$\Is_V (x)$ over all handlebodies $V$ is $\preceq x^{3}$.
\end{cor}

The inequality $\Is_k (x) \preceq x^{k+1}$ was proved by Federer--Fleming \cite{FedererFleming} for integral currents, in Euclidean spaces, and it was later extended by S. Wenger to complete metric spaces with a cone-type inequality \cite{Wenger:isoperim}. For Lipschitz fillings we refer to \cite{ECHLPT}. For fillings of Riemannian hypersurfaces in Banach spaces,  $\Is_k (x) \preceq x^{k+1}$ was proved by Gromov \cite{Gromov(1983)}. We were informed by S. Wenger that the simplicial version of the inequality can also be deduced from \cite{Wenger:isoperim} and \cite[Theorem 1, p.  435]{White}. 

\medskip

The reduction to the subset of round spheres also applies to another type of filling function, the divergence. We briefly recall the definition of divergence, and refer to \cite{BD:HigherFilling1} for more details.

We fix a constant $0<\delta<1$ and an integer $2\leq k\leq n-1$, where $n$ is the dimension of the simplicial complex $X$ we work in. Given a vertex $c$ in $X$, a $k$--dimensional hypersurface $\fh \colon \calm
\to X$ modelled on $\partial V$ such that $k\leq n-1$, and a number $r>0$ that is at most the
distance from $c$ to $\fh (\calm^{(1)})$, \emph{the 
divergence} of this quadruple, denoted $\dv (\fh, c ;r, \delta )$, is 
the infimum of all volumes of
domains modelled on $V$ filling $\fh$ and disjoint from $B(c, \delta r)$.  If
no such domain exists then we set $\dv(\fh, c ;r, \delta )=
\infty$.

\begin{defn}\label{def:kdiv}
The \emph{divergence function modelled on $V$} of the complex $X$, denoted $\Dv_V (r,
\delta )$, is the supremum of all finite values of 
$\dv (\fh, c; r,
\delta )$, where $\fh$ is a hypersurface modelled on $\partial V$ with the distance from $c$ to $\fh (\calm^{(1)})$ at least
$r$ and $\vol (\fh )$ at most $A r^{k}$.

When $V$ is the $(k+1)$--dimensional unit ball, $\Dv_V (r,
\delta )$ is denoted $\Dv^{(k)} (r,
\delta )$, and it is called the $k$--\emph{dimensional divergence function} (or the \emph{$k$--th divergence function}) of $X$.
\end{defn}

In the definition of divergence, as for the isoperimetric function, we fix the constant $A>0$ once and for all, and we do not mention it anymore.

\begin{prop}[\cite{BD:HigherFilling1}]\label{lem:subEucl}
Let $\varepsilon$ and $\delta$ be small enough positive constants.
Assume that $\Is_V (x) \leq \varepsilon x^{k+1}$.
Then $\Dv_V (x , \delta ) = \Is_V (x)$ for every $x$ large enough.	
\end{prop}

Arguments similar to those used for the usual isoperimetry allow to reduce the problem of estimating the divergence to spheres that are round, when a bounded combing exists.

\begin{thm}[\cite{BD:HigherFilling1}]\label{thm:DivRound}
Assume that $X$ is a simplicial complex of dimension $n$ endowed with a bounded quasi-geodesic combing. Let $V$ be a $(k+1)$--dimensional
connected compact sub-manifold of $\R^{k+1}\, $ with connected boundary, where $2\leq k\leq n-1$.

For every $\varepsilon >0$ there exists $\eta >0$ such that the following holds. 

Consider the restricted
divergence function $\Dv^r_V(x, \delta )$, obtained by taking the supremum only over hypersurfaces modelled on $\partial V$ that are
$\eta$--round, of volume at most $2A x^{k}$ and situated outside balls of radius $x$. 

Assume that $\Dv^r_V(x, \delta ) \leq Br^\beta$ for some $\beta \geq k+1$ and $B>0$ universal
constant.  Then the general divergence function $\Dv_V (x, \delta (1- \varepsilon ) )$ is at most $B'r^\beta$ for some $B'>0$ depending on $B, \varepsilon , \eta $ and $X$.
\end{thm}
\subsection{Mapping class groups and  marking complexes}\label{sec:marking}

Let $S$ be a compact oriented surface of genus $g$, with $p$ boundary components, and let $\MCG (S)$ be the mapping class group of $S$. We will use a quasi-isometric model of $\MCG (S)$, the
\emph{marking complex}, $\MM(S)$, defined as follows (see
\cite{MasurMinsky:complex2} for details).  Its vertices,
called \emph{markings}, consist of the following pair of data:
\begin{itemize}
    \item \emph{base curves}: a multicurve consisting of $3g+p-3$
    components, i.e. a maximal simplex in $\CC(S)$.  This collection
    is denoted $\base(\mu)$.

    \item \emph{transversal curves}: to each curve
    $\gamma\in\base(\mu)$ is associated an essential curve.  Letting
    $T$ denote the complexity $1$ component of $S\setminus
    \bigcup_{\alpha\in\base{\mu}, \alpha \neq \gamma}\alpha$, the
    transversal curve to $\gamma$ is a curve $t(\gamma)\in\CC(T)$ with
    $\dist_{\CC(T)}(\gamma,t(\gamma))=1$.
\end{itemize}

Two vertices $\mu,\nu\in\MM(S)$ are connected by an edge if either of
the two conditions hold:

\begin{enumerate}
    \item \emph{Twist}: $\mu$ and $\nu$ differ by a Dehn twist along
    one of the base curves: $\base(\mu)=\base(\nu)$ and all their
    transversal curves agree except for $t_{\mu}(\gamma)$, obtained
    from $t_{\nu}(\gamma)$ by twisting once about the curve $\gamma$.

    \item \emph{Flip}: The base curves and transversal curves of
    $\mu$ and $\nu$ agree except for one pair $(\gamma,
    t(\gamma))\in\mu$ for which the corresponding pair consists of the
    same pair but with the roles of base and transversal reversed.
\end{enumerate}

  Note that after performing one Flip the new base curve may
  intersect several transversal curves. Nevertheless by \cite[Lemma
  2.4]{MasurMinsky:complex2}, there is a finite set of natural ways to
  resolve this issue, yielding a uniformly bounded number of flip
  moves which can be obtained by flipping the pair $(\gamma,
  t(\gamma))\in\mu$; an edge connects each of these possible flips to
  $\mu$.

\begin{thm}[\cite{MasurMinsky:complex2}]\label{MM:marking} The graph
$\MM(S)$ is locally finite and the mapping class group acts
cocompactly and properly discontinuously on it.  In particular, the
orbit map yields a quasi-isometry from $\MCG(S)$ to $\MM(S)$.
\end{thm}

A quasi-geodesic $\fg$ in $\MM(S)$ is
\emph{$\CC(S)$--monotonic} if one can
associate a geodesic $\ft$ in $\CC(S)$ which
{\em shadows} $\fg$ in the sense that if the endpoints of $\fg$ are $\mu$ and $\nu$, then the endpoints of $\ft$ are vertices of
$\pi_{\CC(S)}(\base(\mu))$ and respectively $\pi_{\CC(S)}(\base(\nu))$,
moreover there is a coarsely
monotonic map $v\co\fg\to\ft$ such that $v(\rho)$
is a vertex in $\pi_{\CC(S)}(\base(\rho))$ for every vertex $\rho\in\fg$.

Let $g$ be a pseudo-Anosov on $S$. By
\cite[Proposition~7.6]{MasurMinsky:complex2}, there exists a
quasi-invariant axis of $g$ in $\CC(S)$, that is, a bi-infinite geodesic $\mathfrak a$ in $\CC(S)$ such
that $g^{n} {\mathfrak a}$ is at Hausdorff distance $O(1)$ from $\mathfrak
a$ for all $n\in\Z$.
The set of distances $\dist_{\MM (S)} (\nu , g\nu )$ with
$\pi_{\CC (S)} (\nu )$ at distance $O(1)$ from $\mathfrak a$ admits a
minimum.  Let $\mu$ be a point such that $\dist_{\MM (S)} (\mu , g\mu )$ is
this minimum, and let $\fh $ be a hierarchy path (i.e. a particular
type of quasi-geodesic in $\MM(S)$ as constructed in \cite{MasurMinsky:complex2}) joining $\mu$ and
$g\mu$, such
that $\fh $ shadows a tight geodesic $\ft$ in $\CC (S)$ at distance
$O(1)$ from $\mathfrak a$.  We call the bi-infinite quasi-geodesic
$\p= \bigcup_{n\in \Z } g^n \fh $ a \emph{quasi-axis} of $g$ in $\MM
(S)$.

\me

It follows from \cite{Behrstock:asymptotic} that there exists a
constant $C>0$ such that for every $\delta >0$ small enough, if the
surface $S$ has complexity at least $2$, then $g^{-n}$ and $g^{n}$ may
be joined in $\MCG (S)$ by a path whose length is of order 
$n^2$ which is disjoint
from the $\delta n$--ball around $1$.  Equivalently, for every $\mu
\in \p$, $g^{-n}\mu $ and $g^{n}\mu $ may be joined in $\MM (S)$ by a
path of length of order $n^2$ and disjoint from the $\delta n$--ball
around $\mu $.

\me

We now recall some terminology from \cite[$\S
2.1.8$]{BKMM(2008)}. Let $\Delta$ be an arbitrary simplex in the curve complex $\CC(S)$. Sometimes its set of vertices is called a \emph{multicurve} on $S$. The \emph{open
subsurface} $\mathrm{open} (\Delta)$ determined by $\Delta$ is the union of all the
components of $S\setminus \Delta$ with complexity at least $1$ and of
all the annuli homotopic to curves in $\Delta$. We call
\emph{components of an open subsurface} the list of all the
components and annuli defined as above by a multicurve $\Delta$.  For
$\Delta = \emptyset$ the whole surface $S$ is the unique
component of $\mathrm{open} (\Delta)$.

In what follows we use the notation $W_{1},\ldots ,W_{k}$, with $0< k \leq \xi(S)$, for the list of components of  $\mathrm{open} (\Delta)$, and for every $i\in \{ 1,2,\dots , k\}$ we denote $S\setminus W_{i}$ by $W_i^c$.

Let $\Delta_i$ denote the boundary $\partial W_i$; note that $\Delta$ is the union of the
$\Delta_i$, after removing any duplicate curves.

Define $\QQ(\Delta)$ to be the set of elements of $\MM(S)$ whose
base curves contain $\Delta$.

\begin{lem}\label{lem:qdelta}
\begin{enumerate}
  \item\label{qdel1} If $\mu$ and $\nu$ are two markings in
  $\QQ(\Delta)$ then there exist hierarchy paths joining them and
  entirely contained in $\QQ(\Delta)$.

\me

  \item\label{qdel2} The set $\QQ(\Delta)$
is quasi-isometric to $\MM (W_1)\times \cdots \times \MM(W_k)$.
\end{enumerate}
\end{lem}

\proof  \eqref{qdel1} follows from the construction of hierarchy paths \cite{MasurMinsky:complex2}, while \eqref{qdel2} follows from
\cite[Lemma 2.1]{BehrstockMinsky:rankconj}.\endproof

Given subsets $A_i \subseteq \MM (W_i)$, for $i\in \{1,2,\dots ,k \},$ the
product $A_1\times \cdots \times A_k$ can be identified with a subset
of $\MM (W_1)\times \cdots \times \MM(W_k)$, and by Lemma
\ref{lem:qdelta}, \eqref{qdel2}, it can also be (quasi-)identified
with a subset of $\QQ (\Delta )$.  In what follows we use the same
notation $A_1\times \cdots \times A_k$ for the subsets in $\MM
(W_1)\times \cdots \times \MM(W_k)$ and in~$\QQ (\Delta_i)$.

\me

Let $Z$ be a subsurface of $S$; throughout, all
subsurfaces we consider are implicitly assumed to be essential
subsurfaces. Following \cite{Behrstock:asymptotic}, we consider a
projection $\pi_{\MM(Z)}\co\MM(S)\to 2^{\MM(Z)}$ defined as follows.
For an arbitrary $\mu\in\MM(S)$
we can build a marking on $Z$ by first choosing an element
$\gamma_{1}\in\pi_{\CC (Z)}(\mu)$, and then recursively choosing $\gamma_{n}$
from $\pi_{\CC (Z\setminus \cup_{i<n}\gamma_{i})}(\mu)$, for each $n\leq
\xi(Z)$.  Take these $\gamma_{i}$ to be the base curves of a
marking on $Z$.  For each $\gamma_{i}$ we define its transversal
$t(\gamma_{i})$ to be an element of $\pi_{\mathcal{A}(\gamma_{i})}(\mu)$, where $\mathcal{A}(\gamma_{i})$ is the annulus with core curve $\gamma_i$.  This
process yields a marking, see \cite{Behrstock:asymptotic} for
details. Arbitrary choices were made in this construction, but two
choices in building $\pi_{\MM(Z)}(\mu)$ lead to elements of $\MM(Z)$
whose distance is $O(1)$, where the bound depends
only on $\xi(S)$ \cite{Behrstock:asymptotic}.

\me

Given a marking $\mu$ and a
multicurve $\Delta$, the projection $\pi_{\MM(S\setminus\Delta)}(\mu)$
can be defined as above.  This allows one to construct a
point $\mu'\in \QQ(\Delta)$ which, up to a uniformly bounded error,
is closest to $\mu$.  See
\cite{BehrstockMinsky:rankconj} for details.  The marking $\mu'$ is
obtained by taking
the union of the (possibly partial collection of) base curves $\Delta$
with transversal curves given by $\pi_{\Delta}(\mu)$ together with
the base curves and transversals given by
$\pi_{\MM(S\setminus\Delta)}(\mu)$.  Note that the construction of
$\mu'$ requires, for each subsurface $W$ determined by the multicurve
$\Delta$, the construction of a projection $\pi_{\MM(W)}(\mu )$.  As
explained previously, each $\pi_{\MM(W)}(\mu )$ is
determined up to uniformly bounded distance in $\MM (W)$, thus $\mu'$
is well defined up to uniformly bounded distance, depending only on
the topological type of $S$.

\me

The following is a corollary of the distance formula in
\cite{MasurMinsky:complex2}.

\begin{cor}\label{distsubsurf}
There exist $A\geq 1$ and $B\geq 0$ depending only on $S$ such that
for any subsurface $Z\subset S$, the projection of $\MM (S)$ onto $\MM
(Z)$ is an $(A,B)$--quasi-Lipschitz map, that is for any two markings
$\mu ,\nu \in \MM (S)$ the following holds:
$$
\dist_{\MM (Z)} \left( \pi_{\MM(Z)} (\mu ) \, ,\, \pi_{\MM(Z)} (\nu )
\right) \leq A \dist_{\MM (S)} (\mu , \nu ) +B \, .
$$

Consequently the nearest point projection onto $\QQ (\Delta )$ is a quasi-Lipschitz map.
\end{cor}

Let $g_{i}$ be an element in
$\MCG (S)$ that is pseudo-Anosov when restricted to $W_{i}$ and is the
identity on $W_i^c$, $i\in \{ 1,2,\dots ,k\}$ (this includes the case
of pseudo-Anosov, where $\Delta = \emptyset$). Let $\p_i$ be a quasi-axis of $g_i$ in $\MM (W_i)$.

\begin{prop}\label{prop:pA}
There exists a quasi-Lipschitz map $\Phi_{i}\colon \MM (S)\to \p_{i}$
with the following properties:

\begin{enumerate}
  \item $\Phi_i$ is \emph{coarsely locally-constant} in the complement of
 $\p_{i}\times\MM (W_{i}^c)$, i.e. there exists constants $\lambda >0$ and
$r_0 >0$ with the property that for any point $\mu$ at distance $r
\geq r_0$ from $\p_{i}\times\MM (W_{i}^c)$, the diameter of $\Phi_i
(B(\mu, \lambda r))$ is at most a uniform constant, $c$, which depends
only on $g_i$;

  \medskip

  \item $\Phi_i$ restricted to $\p_{i}\times\MM (W_{i}^c)$ is at uniformly bounded distance from the
projection onto the first component.
\end{enumerate}
\end{prop}

\proof Follows immediately from the proofs of Theorems 3.1 and 3.5 in \cite{BehrstockMinsky:rankconj}.\endproof

The map $\Phi\co
\MM(S)\to \p_{1}\times \cdots \times \p_{k}$ defined by $\Phi (\mu
)=(\Phi_{1}(\mu )\, ,\, \dots \, ,\, \Phi_{k}(\mu ))$ is also quasi-Lipschitz.

The inclusion $\<g_1 \> \times \cdots \<g_k \> \to \MCG (S) $ is a quasi-isometric embedding \cite{FarbLubotzkyMinsky}, and one can use the orbit map to construct a  quasi-isometric embedding $\p_1 \times \cdots \times \p_k \to \MM (S)$, with constants depending on the chosen elements $g_1,\dots , g_k$. The composition of the latter inclusion with the map $\Phi$ is at uniformly bounded distance from the identity, since $\Phi$ is a sort of (quasi-)nearest point projection.

\bigskip

\subsection{Mapping class groups and simplicial complexes}\label{sec:mcgSimplex}

The mapping class group itself and finite index subgroups of it act properly discontinuously cocompactly on various other CW--complexes. Indeed, all finite index torsion-free subgroups of a mapping class group have
classifying spaces given by finite CW-com\-plex\-es, see \cite{Ivanov(1991),Ivanov:mcg}.  For the mapping
class groups themselves, concrete constructions of CW-complexes on
which they act properly discontinuously and with compact quotient
(complexes that are moreover cocompact models for classifying spaces
for proper actions) are described in \cite{JiWolpert,
Mislin:classifying}.  Any of these CW-complexes can be
used to define filling functions for the mapping class group.

Another approach is to apply Theorem 10.2.6 in \cite{ECHLPT}.  The
mapping class groups are automatic \cite{Mosher:automatic}, hence
combable, so they are $\mathcal{F}_\infty$.  In particular, for every
$n\geq 0$, they act by simplicial isomorphisms, properly
discontinuously, with trivial stabilizers of vertices, on a simplicial complex $X$ of dimension $n+1$ and
$n$-connected, such that the quotient has finitely many cells.

For every surface $S$ we denote by $X_S$ a simplicial complex with
properties as above.  Note that we allow the case when $S$ has several
connected components.

Using the previous section and \cite[Lemma 12]{AWP}, the following
properties can be established about these simplicial complexes, up to
repeated barycentric subdivisions.  We use the same generic objects
defined in the previous section, with the same notation and
terminology.

Given a pseudo-Anosov $g\in\MCG (S)$, there exists a bi-infinite
almost geodesic $\widehat{\p }$ in $X_S^{(1)}$ such that $g$ acts on
it with compact quotient, and the points on it quasi-minimize the
displacement by $g$.  Such a path is a \emph{quasi-axis} of $g$.
Moreover, for every $x \in \widehat{\p}$, $g^{-n}x $ and $g^{n}x $ may
be joined in $X_S$ by a curve of length of order $n^2$ which is disjoint from
the $\delta n$--ball around $x$.

For $\Delta, W_1, \dots, W_k$ defined as previously, there exists a
quasi-isometric embedding which is also a simplicial map
$\Lambda_\Delta \co X_{W_1} \times \cdots \times X_{W_k} \to X_S$, with
image a subcomplex such that every two points in it can be joined in
it by a path that is an almost geodesic in $X_S$.

As before, given subsets $B_i$ in $X_{W_i}$, $i=1,2,\dots, k\, ,$ we
let $B_1 \times \cdots \times B_k$ denote the subset in $X_{W_1}
\times \cdots \times X_{W_k}$ and its image in $X_S$.

There exists a bounded perturbation of the nearest point projection
$\widehat{\pi}_\Delta \co X_S \to X_{W_1} \times \cdots \times X_{W_k}$
which is a simplicial map, and a quasi-Lipschitz map.
This allows to define a map with the same properties
$\widehat{\pi}_{W_i} \co X_S \to X_{W_i}$ for $i=1,\dots, k\, $.
Consider pure reducible elements $g_i$ as before, and their quasi-axes
$\widehat{\p}_i $ in $X_{W_i}$.

The projections $\Phi_i$ defined in Proposition~\ref{prop:pA} allow to define
a simplicial map $\widehat{\Phi}_i \co X_S \to \R $ that is
quasi-Lipschitz and coarsely locally-constant in the complement of
$\widehat{\p}_i \times X_{W_i^c}$, while its restriction to
$\widehat{\p}_{i}\times\MM (W_{i}^c)$ is at uniformly bounded distance
from the projection onto the first component.
We can then define the map $\widehat{\Phi} \co X_S \to \R^k $,
$\widehat{\Phi} = \left( \widehat{\Phi}_1, \dots, \widehat{\Phi}_k
\right)$ which is also quasi-Lipschitz.

The almost geodesics $\widehat{\pgot}_i \co \R \to X_{W_i}$ (which can
also be seen as simplicial maps) define an inclusion
$\Upsilon_\Delta\co
\R^k \to X_{W_1} \times \cdots \times X_{W_k}$ which is simplicial,
a quasi-isometric embedding, and equivariant with respect to the action of
$\langle g_1 \rangle \times \langle g_2 \rangle \times \cdots \times
\langle g_k \rangle $. Hence, we have an inclusion $\Upsilon \co \R^k \to X_S$
,$\Upsilon = \Lambda_\Delta \circ \Upsilon_\Delta$ with the same
properties; moreover, we also have that  $\widehat{\Phi} \circ \Upsilon$ is at
uniformly bounded distance from the identity.

\section{Higher dimensional filling and divergence in mapping class groups}\label{sec:divmcg}

In this section we fix an arbitrary surface $S$ with $\xi(S)>0$. As in Section \ref{sec:mcgSimplex}, we consider a simplicial complex $X_S$ on which $\MCGS$ acts properly discontinuously and such that the quotient of each skeleton $X_S^{(m)}$ is composed of finitely many $m$--simplices. From  the fact that mapping class groups are automatic
\cite{Mosher:automatic} it follows that $X_S$ has a bounded $(L,C)$--quasi-geodesic combing.

\subsection{Isoperimetry in mapping class groups}\label{sec:fillingmcg}

Known facts about mapping class groups imply that
 for each dimension $k$ below the quasi-flat rank, the isoperimetric
 function in $X_{S}$ is asymptotically equal to the
 \emph{$k$--dimensional Euclidean isoperimetric function}, and that above
 the rank the isoperimetric function is sub-Euclidean.

\begin{thm}\label{thm:genus3}
The $k$--th isoperimetric function in the mapping class group of a surface satisfies
$\Is_k (x) \asymp x^{k+1}$ for $k< \xi (S)$ and
$\Is_k (x)=o( x^{k+1})$ for $k\geq \xi (S)$.
\end{thm}

\proof
By Theorem \ref{cor:rounddec}, for every integer $k\geq 1$ we have
$\Is_{k} (x) \preceq
x^{k+1}$.  Moreover, the presence of
quasi-flats of dimension $k$ inside the mapping class groups (see
Proposition \ref{prop:cube} for a construction of some such
quasi-flats and Corollary~\ref{distsubsurf} for relevant results about
their geometry) and \cite[Theorem 2]{AWP} imply that for $k<\xi (S)$ we have
$\Is_k (x) \asymp x^{k+1}$.

For $k\geq \xi(S)$, the Theorem follows from \cite[Theorem 1, p.  435]{White}, \cite{Wenger:isoperim}, and the fact that the maximal dimension of locally
compact subsets in an arbitrary asymptotic cone of $\MCG (S)$ is $\xi
(S)$ \cite{BehrstockMinsky:rankconj}.

For the convenience of the reader, we also provide a more explicit argument.
According to \cite{BD:HigherFilling1}, it suffices to
prove that the filling function is $o(x^{k+1})$ only for spheres (respectively surfaces) that are round and unfolded.  We argue for a contradiction and assume
that there exists a sequence $\rgot_n$ of $k$-dimensional spheres for $k\geq 3$ (respectively of surfaces for $k=2$) that
are round and unfolded, of volume $\asymp x_n^{k}$ and of filling volume at
least $\lambda x_n^{k+1}$, where $\lambda$ is a positive constant and
$x_n \to \infty\, $.  Let $\dgot_n$ be filling $(k+1)$-dimensional
balls (respectively filling handlebodies) realizing $\fvol (\rgot_n )$ and with a minimal number of chambers in the domain.

The argument in \cite[pp. 263--264]{Wenger:asrk} with
$T_n= \rgot_n$ and $S_n = \dgot_n$ implies that the sequence $\left( \dgot_n \right)$  yields a compact subset of dimension $k+1$ in an asymptotic cone of $X_S$, a contradiction.\endproof

For  $k\geq \xi (S)$ we conjecture $\Is_k (x) \asymp x^{k}$. In a forthcoming paper \cite{BehrstockDrutu:Linear} we prove that an asymptotic version of this holds. 
The sharp 
result holds in low genus:

\begin{thm}\label{thm:genus012}
Given a surface $S$ of genus $0$ or $1$, or of genus $2$ and without
boundary, the $k$--th isoperimetric function in the mapping class group of
$S$ satisfies $\Is_k (x) \asymp x^{k}$ for $k\geq \xi (S)$.
\end{thm}

\proof It was recently established that the mapping class group of a surface has a cocompact classifying space
for proper actions of dimension equal to the virtual cohomological
dimension \cite{HenselOPr,AramayonaPerez}.  The virtual cohomological
dimension for the mapping class group of $S$ is: $p-3$ if $g=0$;
$4g-5$ if $p=0$; and $4g+p-4$ if both $g$ and $p$ are positive
\cite{Harer:vcd}.

Since the surfaces in the hypothesis of the theorem thus have
$vcd(\MCG(S))=\xi(S)$, the result then follows from the fact that
if a group has a cocompact classifying space for proper actions
of dimension $r$, then $\Is_k (x) \asymp x^{k}$ for all $k\geq r$
\cite[Corollary 9]{AWP}.
\endproof

\medskip

\subsection{Divergence in  mapping class groups}

In the mapping class group, the value for rank analogous to that in a
symmetric space is the quasi-flat rank, i.e. the maximal dimension of
a quasi-flat in the Cayley graph of the group.  As discussed
previously, for $\MCG(S)$ this
rank is given by $\xi(S)=3g+p-3$ \cite{BehrstockMinsky:rankconj,
Hamenstadt:qirigidity}.  Below we show how, analogous to the case of
symmetric spaces, this rank plays a critical
role for divergence in mapping class groups as well.

\begin{thm}\label{thm:divOverRk}
Given a surface $S$  and an arbitrary integer
$k \geq \xi(S)$, the $k$-dimensional divergence in $\MCG(S)$
satisfies $\Dv_{k} (x)= o(x^{k+1})$.

If moreover $S$ is of genus $0$ or $1$, or of genus $2$ and without
boundary, then $\Dv_{k} (x) \asymp x^{k}$.
\end{thm}

\proof The result follows from Proposition \ref{lem:subEucl}, and from
Theorems \ref{thm:genus3} and \ref{thm:genus012}.  \endproof

\begin{thm}\label{thm:divsubrank} For any $S$ and for any integer
$0\leq k < \xi(S)$, the $k$-dimensional divergence in $\MCG(S)$
satisfies $\Dv_{k}\succeq x^{k+2}$.
\end{thm}

The rest of the section is devoted to the proof of Theorem \ref{thm:divsubrank}. We use the terminology introduced in Sections \ref{sec:marking} and
\ref{sec:mcgSimplex}.  In particular we consider a compact connected
orientable surface $S$ of complexity $m=\xi (S) \geq 2$, and an
 $m$-connected simplicial complex $X_S$ of dimension $m + 1$
on which $\MCG (S)$ acts by simplicial isomorphisms,
properly discontinuously, with trivial stabilizers of vertices and
such that the quotient has finitely many simplices.
Theorem~\ref{thm:divsubrank} is a direct consequence of
Proposition~\ref{prop:cube} below. Indeed, in this proposition we show that for every $0< k
\leq \xi(S)$ there exist naturally arising $(k-1)$--dimensional spheres
in $X_S$ which have divergence $\asymp x^{k+1}$.

\begin{prop}\label{prop:cube}
Let $\Delta$ be a multicurve on $S$ and let $W_{1},\ldots ,W_{k}$ be the components
of $\mathrm{open} (\Delta)$, $0< k \leq \xi(S)$. For $i=1,\dots , k$, consider $g_{i}$ elements in
$\MCG (S)$ that are pseudo-Anosov when restricted to $W_{i}$ and the
identity on $W_i^c=S\setminus W_{i}$. Let $\Upsilon
\co\R^k \to X_S$ be a simplicial map defining a quasi-flat and equivariant for the action of $\langle g_1 \rangle \times \langle g_2 \rangle \times
\cdots \times \langle g_k \rangle $, and let $x_0 = \Upsilon (0)$.

Let $\cals_n^{(k-1)}$ denote the boundary of the minimal simplicial subcomplex of $\R^k$ covering the cube $\{ (x_1, \dots , x_k)\in \R^k\; ;\; |x_i| \leq
n, i=1,2,...k \}$ and let $\sgot_n^{(k-1)}$ be the $(k-1)$--dimensional sphere in $X_S$ defined by the restriction of $\Upsilon$ to $\cals_n^{(k-1)}$.

\begin{enumerate}
  \item\label{least} There exists a constant $C_1>0$, such that
  for every small enough $\delta \in (0, 1)$, every $k$--disk
  $\dgot$ filling the sphere $\sgot_n^{(k-1)}$ in $X_S$ and
  disjoint from the $\delta n$--ball around $x_0$ has area at least $C_1
  n^{k+1}$.

\medskip

  \item\label{most} There exists a constant $C_2>0$ such that, if
  at least one of the surfaces $W_i$ has
  complexity $\geq 2$ (which implies, in particular, that $k<
  3g+p-3$),
  then for every $\delta >0$ small enough there exists a $k$--disk
  $\dgot$ filling the sphere $\sgot_n^{(k-1)}$, disjoint
  from the $\delta n$--ball around $x_0$, and with area at most $C_2 n^{k+1}$.
\end{enumerate}
\end{prop}

\begin{rmk}
Proposition \ref{prop:cube}, (\ref{least}), implies that the
$k$--dimensional quasi-flats $\Upsilon (\R^k)$ in $X_S$ and $\langle
g_1 \rangle \times \langle g_2 \rangle \times \cdots \times \langle
g_k \rangle $ in $\MCG(S)$ are maximal in the sense that no quasi-flat
of strictly larger dimension quasi-contains them.  This is a subtle
point, even though it is
obvious when considering quasi-flats of the same type.
\end{rmk}

\proof  (\ref{least}) Let $\dgot :\dd \to X_S$ be a $k$--dimensional disk which fills
$\sgot_n^{(k-1)}$ in $X_S$ and is disjoint from the $\delta n$--ball around
$x_0$. Then $\widehat{\Phi} \circ \dgot $ is a disk filling $\widehat{\Phi} \circ \sgot_n^{(k-1)}\, $
in $\R^k$. Since $\widehat{\Phi} \circ \Upsilon$ is at uniformly bounded distance from the identity, $\widehat{\Phi} \circ \sgot_n^{(k-1)}$ is at uniformly bounded distance from $\cals_n^{(k-1)}$.

In particular, for $\epsilon >0$ small  enough, the image of $\widehat{\Phi} \circ \dgot$ covers
the full cube in $\R^k$ centered in $x_0$ and of edge length
$\epsilon n$. We denote the latter by $\mathrm{Cube}_\epsilon$ and the initial full cube of edge length $2n$ by $\mathrm{Cube}_2$.

For ease of notation, we explain the argument in the case $k=2$ and
then, after each step, how it can be modified to yield the general case, when the required modifications are not obvious.

For any choice of $\eta<n/2$ we may consider a grid subdividing $\mathrm{Cube}_2$ into squares with edge length $2\eta$ (to simplify the
discussion we will assume that both $n$ and $\epsilon n$ are integral multiples of $2\eta$;
otherwise an additional discussion of boundary rectangles is needed,
which adds only notational complications, see, e.g. \cite{LangSchroeder})
and for each of the squares
composing it consider the sub-square of edge $\eta$ obtained by
shrinking with factor $\frac 12$ around the center.  Let $\Pi$ be one
such full square that is moreover contained in $\mathrm{Cube}_\epsilon$.  Consider the $\ell_1$--geodesic $\pg$ starting in the
midpoint $a$ of the upper horizontal edge of the square, going through
the center of the square and ending in the midpoint $b$ of the right
hand side vertical edge.

In what follows we denote the map $\widehat{\Phi} \circ \dgot : \dd \to \R^k$ by $f$.

\me

{\bf Step 1.} The first step is to prove that for some constant
$\lambda\ll \eta$ the set $f\iv (\pg )$ has a connected component $\mathfrak C$ for
which $f({\mathfrak C})$ intersects the $\lambda$--ball around each of the two points $a$ and
$b$ respectively. This argument is
inspired from the proof of \cite[Proposition 3.2]{LangSchroeder}.

Assume on the contrary that this is not true.  Since the map $f$ is
simplicial, defined on a triangulation of $\mathbb{D}^2$, there are finitely
many connected components of $f\iv (\pg )$.  According to our
hypothesis all the images of components that intersect $B(a ,\lambda
)$ do not intersect $B(b, \lambda )$.  Let $G$ be the union of the
components with image by $f$ intersecting $B(a ,\lambda )$ but not
$B(b ,\lambda )$.  Since both $G$ and $f\iv (\pg) \setminus G$ are
compact there exists $\varepsilon >0$ such that $U_\varepsilon = \{
x\in \mathbb{D}^2\; ;\; \dist (x, G ) < \varepsilon \}$ intersects
$f\iv (\pg )$ only in $G$, and $f (U_\varepsilon )$ intersects the
boundary of the square $\Pi $ only in the upper horizontal edge of the
square.

In what follows we construct a continuous map $F$ from the 2-dimensional disk $\mathbb{D}^2$ to $\R^2$, which coincides with $f$
outside $U_\varepsilon$ and on $f\iv (\partial \Pi )$ (here $\partial \Pi$ is the boundary of $\Pi$) and with image
not containing the geodesic segment ${\mathfrak{c}}=\left[ \pg
\setminus \{ a\} \right] \cap B(a, \lambda )\, $.  This will finish
this step, since it will imply that $\partial \Pi$, contractible in $F(\mathbb{D}^2)$, is therefore contractible in $\R^2 \setminus {\mathfrak{c}}$, a contradiction.

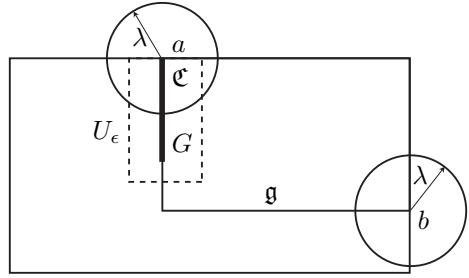
\begin{figure}[h]
\includegraphics[width=0.5\textwidth]{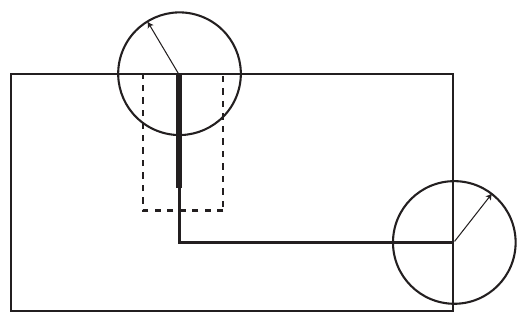}
\put  (-24, 38) {$\lambda$}
\put  (-130, 90)  {$\lambda$}
\put  (-115, 88) {$a$}
\put   (-22, 21){$b$}
\put  (-115, 75){${\mathfrak C}$}
\put   (-115, 50){$G$}
 \put  (-145, 55){$U_{\epsilon}$}
\put  (-80, 32){${\mathfrak g}$}
\caption{The set-up for Step 1}\label{fig:step_1}
\end{figure}

Consider the map $F \colon \mathbb{D}^2 \to \R^2$ defined as follows.
For $x\not\in U_\varepsilon$ let $F(x)=f(x)\, $. For $x$ in
$U_\varepsilon$ let $\pg_x$ be the $\ell_1$ geodesic through $f(x)$, contained in $\Pi$, composed of a vertical and a horizontal segment with common
endpoint on the diagonal from the upper right corner to the lower left
corner. Let $a_x$ be the intersection point of $\pg_x$ with the upper
horizontal edge of $\partial \Pi$, and let $d(x)$ be the $\ell^1$--distance from $a_x$ to $f(x)$.

Let $\ell$ be the length of the geodesic $\pg $ and let $t(x)=
\ell \left[ 1-\frac{\dist (x, G )}{\varepsilon }\right]$.  We then
define $F(x)$ to be the point on $\pg_x$ with distance to the upper
horizontal edge along $\pg_x$ equal to the maximum between $0$ and $d(x)-t(x)\, $.

In other words, if $d(x)$ is at most $t(x)$, then $F(x)$ equals $a_x$; while if $d(x)$ is larger than $t(x)$, then $F(x)$ equals the point on $\pg_x$ at distance $d(x)-t(x)$ from $a_x$.

The map $F$ is continuous, coincides with $f$ outside $U_\varepsilon$ and on the
boundary of $\Pi$, and its image does not contain the piece of
geodesic ${\mathfrak{c}}=\left[ \pg \setminus \{ a\} \right] \cap B(a,
\lambda )$. This provides a contradiction which, as noted above, finishes the step when
$k=2$.

In the general case of dimension $k$ one may consider, instead of the
upper right square bounded by $\pg$ and the boundary of $\Pi$, a cube
of edge length half the length of the edge of $\Pi \, $ with one
vertex the center of $\Pi$ and half of its faces contained
in $\partial \Pi$.  The geodesic $\pg$ is to be replaced by the
$\ell_1$ geodesic joining the endpoints of a big diagonal not containing
the center; denote this new geodesic also by $\pg$.  It suffices to
prove that $f\iv (\pg)$ has one connected component $\mathfrak C$ such that
$f({\mathfrak C})$ intersects the $\lambda$ neighborhoods of both endpoints of
$\pg$.  This is completed as above by showing that for each cube $\Pi$
there exists a connected set ${\mathfrak C}_\Pi$ in $f\iv (\pg )$ such that
$f\left( {\mathfrak C}_\Pi \right)$ covers the $\ell^1$ geodesic constructed as
above, except eventually the $\lambda$--neighborhoods of its endpoints.

\me

{\bf Step 2.} We denote by $\pg_\Pi$ the subgeodesic of $\pg$ covered by $f({\mathfrak C}_\Pi)$.  Let
$V_\Pi$ be the set of vertices of simplices in the triangulation $\dd$ of
$\mathbb{D}^{2}$ that intersect ${\mathfrak C}_\Pi$.  Since the image by $f$ of
every simplex has diameter at most $1$ it follows that the Hausdorff
distance between $f(V_\Pi )$ and $\pg_\Pi$ is at most $1$.

Note that since $C_\Pi$ is a connected set, the set $\dgot (V_\Pi)$ is
$1$--coarsely connected in the sense that: for every two vertices
$x,y\in \dgot (V_\Pi)$, there exists a finite sequence of vertices
$z_1=x,z_2,...z_n=y$ such that for each $i$ we have $z_{i}\in \dgot
(V_\Pi)$ and $\dist (z_i, z_{i+1})\leq 1$.

In each $V_\Pi$ we fix a vertex $x_\Pi$ such that $\widehat{\Phi} \circ \dgot
(x_\Pi )$ is at distance at most $1$ from the center of $\Pi$.  Let
$\P$ be the set of all shrunken squares $\Pi$ that appear in the grid
and are contained in $\mathrm{Cube}_\epsilon$, and let $\V $ be the set
of vertices $\{ x_\Pi \; ;\; \Pi \in \P \} \, $.  Note that there are
approximately $\frac{\epsilon^2 n^2}{\eta^2}$ elements in $\V$.  Fix a small constant $\beta>0$.

\me

\noindent \emph{Case 1.}\quad Assume that at least half of the points
in $\V $ have the property that their images by $\dgot$ are at distance
$\geq \frac{\beta}{2} n$ from $X_{W_1}\times \cdot \times X_{W_k}$.  Denote this subset of
$\V $ by $\V'$.

\me

If a vertex $w$ is at distance $\geq \frac{\beta}{2} n$ from $X_{W_1}\times \cdot \times X_{W_k}$ then there exists $i\in \{ 1,2,...k \}$ such that $w$ is at
distance $\geq \alpha n$ from $X_{W_i}\times X_{W_i^c}$, where $\alpha$ is a
fixed constant depending only on $\beta$ and $\xi(S)$.  This follows
from the corresponding statement in the marking complex: if a point $\mu$ is at distance $\geq \frac{\beta'}{2} n$ from $\QQ
(\Delta )$ then there exists $i\in \{ 1,2,...k \}$ such that it is at
distance $\geq \alpha' n$ from $\QQ (\Delta_i )$, where $\alpha'$ is a
fixed constant depending only on $\beta'$ and $\xi(S)$.  The latter statement can be proved by induction on $\xi(S)$ and the fact that if $\mu$ is at
distance $\leq \alpha n$ from all $\QQ (\Delta_i )$ then its
projection on $\QQ (\Delta_1)$ is at distance $\leq L \alpha n$ from
$\QQ (\Delta_1 \cup \Delta_i )$ and hence $\mu $ is at distance $\leq
(L+1) \alpha n$ from $\QQ (\Delta_1 \cup \Delta_i )$.

\me

We may thus assume that at least $\frac{1}{k}$ of the points in $\V'$
are such that their image by $\dgot$ is at distance $\geq \alpha n$ from
$X_{W_i}\times X_{W_i^c}$ for a fixed $i$.  Note that for each of these
vertices $x_\Pi$ the corresponding set $\dgot (V_\Pi)$ contains a point
such that its image by $\widehat{\Phi}_i$ is at distance at least $\eta$ from
the $\widehat{\Phi}_i$--image of $x_\Pi\, $.  This, the fact that $\widehat{\Phi}_i$ is
coarsely locally-constant (see Proposition~\ref{prop:pA}) and the fact that $x_\Pi$ is at distance at
least $\alpha n$ from $X_{W_i}\times X_{W_i^c}$ implies that $\dgot (V_\Pi)$
contains a point at distance at least $\lambda \alpha n$ from $x_\Pi
\, $.  Since $\dgot (V_\Pi)$ is $1$--coarsely connected, it follows that
$\dgot (V_\Pi)$ contains at least $\lambda \alpha n\, $ distinct points.
Note that if $\Pi$ and $\Pi'$ are two disjoint cubes, then for $\eta$
large enough $\dgot (V_\Pi)$ and $\dgot (V_\Pi')$ are disjoint.  We have
thus obtained that the image of $\dgot$ contains at least $\frac{1}{2k}
\frac{\epsilon^2 n^2}{\eta^2 } \cdot \lambda \alpha n d$ distinct points, 
images of interior vertices. We proved in \cite{BD:HigherFilling1} 
that with these hypotheses the chambers are all non-collapsed
and thus the area of
$\dgot$ is $\succeq n^3$.

\me

\noindent \emph{Case 2.}\quad Now assume that at least half of the
points in $\V$ have the property that their $\dgot$--images are at
distance $\leq \frac{\beta}{2} n$ from $X_{W_1}\times \cdots \times X_{W_k}$. Let $\V'$ be this new subset of vertices. Since $\dgot (\V')$ avoids a $\delta n$--ball
around $x_0$ and $\widehat{\Phi } \circ \dgot (\V')$ is contained in
$\mathrm{Cube}_\epsilon$, it follows that, for $\beta$ and $\epsilon$ small enough, by post-composing $\dgot$ with the projection onto $X_{W_1}\times \cdots \times X_{W_k}$, we obtain a filling disk $\dgot' = \widehat{\pi}_\Delta \circ \dgot$ such that the points in $\V' $ are sent at distance $\geq \frac{\beta}{2} n$
from $\Upsilon_\Delta (\R^k )$.  Then at least $\frac{1}{k}$ of the
points in $\V'$ have image by $\dgot'$ at distance $\geq
\frac{\beta}{2k} n$ from $\widehat{\p}_i$ for some fixed $i\in \{ 1,2,...k\}$.
Denote this subset of $\V' $ by $\V''$.  For notational simplicity we
assume that $i=1$.

 For every vertex $x_\Pi$ in $\V''$ there exists a point in $V_\Pi$
 such that given their respective images by $\dgot'$, the projections of
 their respective first coordinates onto $\widehat{\p}_1$ are at distance at
 least $\eta\, $.  Then with the same argument as in Case 1 we deduce
 that $\dgot' (V_\Pi )$ contains at least $\lambda \alpha n$ points, hence
 $\dgot'$ (and therefore $\dgot$, up to multiplication by some universal constant in $(0,1)$) has area at least $\frac{1}{6k}
\frac{ \epsilon^2 n^2}{\eta^2} \cdot \lambda \alpha n \, $. 

\medskip

(\ref{most}) Assume that $W_1$ has complexity $\geq 2\, $.  The point
$x_0$ is the image by $\Lambda_\Delta$ of a point $y_0 = (y_0^1, \dots
, y_0^k)$ in $\Upsilon_\Delta (\R^k )\, $.  According to the arguments
in Section \ref{sec:mcgSimplex}, $\widehat{\pgot}_1 (-n)$ and
$\widehat{\pgot}_1 (-n)$ can be joined in $X_{W_1}$ by a path $\cgot$
with length in $[a n^2 , bn^2]$, where $0<a<b$, path avoiding the ball
centered in $y_0^1$ and of radius $\delta' n\, $.  Let $\dgot$ be a
disk of area whose area is an affine function of $n^3$ 
filling the loop $\cgot \cup
\widehat{\pgot}_1([-n, n])$.  It can be obtained by taking, for every
point $x$ on $\cgot $ its projection $x'$ on $\widehat{\pgot}_1([-n,
n])$, letting $x$ vary, using the fact that the projection on
$\widehat{\pgot}_1$ is quasi-Lipschitz and coarsely locally-constant
(see Proposition~\ref{prop:pA}) and quasi-geodesics $\q_{x'x}$ given
by the combing with basepoint $x'$.

 Then $\cgot \times \widehat{\pgot}_2([-n, n]) \times \cdots \times \widehat{\pgot}_k([-n, n]) \cup \bigcup_{i=2}^k \dgot \times \widehat{\pgot}_2([-n, n]) \times \cdots \times \widehat{\pgot}_{i-1}([-n, n])\times  \widehat{\pgot}_{i}(\pm n) \times \widehat{\pgot}_{i+1}([-n, n])\times \cdots \times \widehat{\pgot}_k([-n, n])$ compose a filling disk for the $\Upsilon_\Delta$--image of the boundary of the cube of edge $2n$ centered in $0$, moreover this disk has area
$n^{k+1}$ and it is disjoint from the $\delta'' n$--ball
around~$y_0$, for an appropriate choice of $\delta''$ and $\delta'$. By applying $\Lambda_\Delta$ we obtain a disk filling the given sphere, disjoint from the $\delta n$--ball
around~$x_0$ and of area
$n^{k+1}$. \endproof

\section{Higher dimensional divergence of $CAT(0)$--groups}
\label{sec:divcat}

In this section, we show how the above technique for computing higher
divergence in mapping class
groups can be applied in the context of $CAT(0)$--groups to obtain
interesting examples. In particular, in this section we
obtain examples that further clarify the connection
between the quasi-flat rank and the divergence rank in
$CAT(0)$--groups, by describing the behavior of the higher dimensional
divergence when the dimension is below the rank.

\begin{thm}\label{thm:divmcg}
For every positive integers $r$ and $n$ there exist finitely generated
$CAT(0)$-groups with quasi-flat rank $2r$ (equal to the flat rank of
the $CAT(0)$--space they act cocompactly on) and such that the
$(r-1)$-dimensional divergence satisfies $\Dv_{r-1}\asymp 
x^{r+n}$.
\end{thm}

\proof A proof of this result in the case $r=1$ can be found in
\cite[Theorem 1.1]{BehrstockDrutu:thick2}.  The study of higher rank divergence
is more delicate; the main additional tool needed to extend from the $r=1$
case to higher dimensions is obtained by adapting 
the arguments we developed above to compute higher dimensional divergence in
the mapping class group.

In \cite[Proposition 5.2]{BehrstockDrutu:thick2} we described an iterative
construction of a compact $CAT(0)$ space $M_{n}$ of rank $2$ which we now recall. Our construction starts with a
3-dimensional non-geometric graph manifold, $M_{1}$, with a periodic Morse
geodesic, $\fg_1$, such that
any of the lifts $\widetilde{\fg}_1$ to the universal cover
$\widetilde{M}_1$ has divergence of order $x^{2}$. It was shown in
\cite{Gersten:divergence3} that
such examples exist. The space $M_{n+1}$ is
obtained by taking two copies of $M_{n}$ and amalgamating along
their corresponding copies of $\fg_n$. The key element in the proof of
\cite[Proposition 5.2]{BehrstockDrutu:thick2} is
that there exists $\eta_n >0 $ such that, given two points
$\alpha$ and $\alpha'$ in the universal cover $\widetilde{M}_n$, at distance at least $x$
from a lift $\widetilde{\fg}_n$ of $\fg_n$, and with projections on $\widetilde{\fg}_n$
at distance at least $\eta_n$, the shortest path ${\mathfrak{c}}'$
joining $\alpha$ and $\alpha'$ outside the $x$--tubular neighborhood
of $\widetilde{\fg}_n$ has length $\succeq x^{n+1}$.

Let $r$ and $n$ be arbitrary fixed integers. We fix a compact $CAT(0)$ space $M_n$,
together with a Morse periodic geodesic $\fg_n$, as above. Take $N_n$ to be a cartesian product of $r$ copies of
$M_n$, which yields a space of rank $2r$.  Let $\fg_n\times \cdots
\times \fg_n$ be the $r$--dimensional torus with factors $r$ copies of
the periodic geodesic above mentioned.  We will
restrict our attention to those triangulations on $M_n$ which extend
triangulations of $\fg_n$ and of the Seifert components of the space
$M_{1}$, and we will consider
triangulations on $N_n$ that refine the product structure on this
space.

The argument in the proof of Theorem \ref{prop:cube},
Part~(1), Step~2, Case~2, carries over almost verbatim to give that, in the universal cover
$\widetilde{N}_n$, the boundary of a cube of edge length 
$x$ in an
$r$-dimensional flat $\widetilde{\fg}_n\times \cdots \times
\widetilde{\fg}_n$, where $\widetilde{\fg}_n$ is a lift of  $\fg_n$, has $(r-1)$--dimensional divergence 
$\asymp x^{r+n}$ with respect to its
center, for any choice of $\delta \in (0,1)$.

We now prove, by induction on $n$, that the $(r-1)$--dimensional
divergence in $\widetilde{N}_n$ is at most $x^{r+n}$. In fact, we prove the stronger statement that for every $r$--dimensional
connected compact sub-manifold $V$ of $\R^{r}\, $, $V$ with connected boundary, the corresponding divergence function $\Dv_V$ is at most $Bx^{r+n}$, where $B$ is independent of $V$.

According to Theorem \ref{thm:DivRound}, it suffices to prove the above statement for hypersurfaces that are $\eta$--round, for some fixed small $\eta$.

As in \cite{BehrstockDrutu:thick2}, we choose $M_1$ to be a
particularly simple graph manifold, namely, we let this be the space
constructed by taking a pair of hyperbolic surfaces each with at least
one boundary component, crossing each with a circle, and gluing the
two $3$--manifolds $M_0$ and $M_0'$, thus obtained, along a boundary
torus by flipping the base and fiber directions.  This implies that
$N_1$ is obtained by gluing the product $N_0$ of $r$ copies of $M_0$
with the product $N_0'$ of $r$ copies of $M_0'$ along the product of
$r$ copies of a $2$--dimensional torus, following a flip in each
factor.
As $M_{n+1}$ is inductively obtained by isometrically gluing two
copies of $M_{n}$ along $\fg_n$, we have that $N_{n+1}$ is obtained from $N_n$ by taking
two copies of it, and identifying the two respective copies of the
$r$--torus $\fg_n\times \cdots \times \fg_n$.

In the initial step of the induction, the two symmetric spaces
$\widetilde{N}_0$ and $\widetilde{N}_0'$ of rank $2r$ have
$(r-1)$-dimensional divergence $\leq B x^{r}$. The same bound holds for every $\Dv_V$ with $V$ of dimension $r$, and the constant $B$ independent of $V$.

One has to deduce
from the initial step that the $(r-1)$--dimensional divergence of
$\widetilde{N}_1$ is $\preceq x^{r+1}$.  In the inductive step, from
the fact that $\widetilde{N}_n$ has $(r-1)$--dimensional divergence
$\preceq x^{r+n}$ it must be deduced that $\widetilde{N}_{n+1}$ has
$(r-1)$--dimensional divergence $\preceq x^{r+n+1}$.  The two
arguments are very similar, so we explain the latter only.

Consider a point $c$ in $\widetilde{N}_{n+1}$ and an
$(r-1)$--dimensional hypersurface $\fh$ outside $B(c,x)$ and of volume
at most $Ax^{r-1}$.  The space $\widetilde{N}_{n+1}$ is composed of
copies of $\widetilde{N}_{n}$ glued along flats of dimension $r$, and
this decomposition is encoded by the Bass--Serre tree of the
fundamental group: vertices correspond to copies of
$\widetilde{N}_{n}$, while edges correspond to gluings.  In what
follows we call the copies of $\widetilde{N}_{n}$ \emph{geometric
components}, or, simply, \emph{components} of $\widetilde{N}_{n+1}$,
following the terminology in $3$-dimensional manifolds.  We call \emph{separating flats} the
$r$-dimensional flats that are lifts of the $r$--torus along which the
gluing is done.  Given a vertex $v$ in the
Bass-Serre tree, we use the notation $N(v)$ to designate the geometric component
corresponding to that vertex.  Likewise, given an edge $e$ we denote
by $F(e)$ the separating flat corresponding to it.

 For an arbitrary convex set $C$ in $\widetilde{N}_{n}$, we denote by
 $\pi_C$ the nearest point projection onto $C$.

 Note that every separating flat splits $\widetilde{N}_{n+1}$ into two
 open convex subsets, and their closures (obtained by adding the flat
 in question) are likewise convex.  We call the latter closures
 \emph{half spaces}.  Given a point $p$ contained in only one
 geometric component, a half space \emph{opposite to it} is one that does
 not contain that point.

 Without loss of generality we may assume that the point $c$ is a
 vertex of the triangulation and it is contained in only one
 component, therefore only one vertex $v_c$ in the Bass-Serre tree
 corresponds to it.  The hypersurface $\fh$ determines a finite
 connected sub-tree $T_\fh$ inside the Bass-Serre tree, the edges of
 which are all the separating flats crossed by $\fh$, and the vertices
 of which are all the components intersecting $\fh$ in an
 $(r-1)$--dimensional hypersurface with boundary.

Given a separating flat $F$, we denote by $\fh_F$ the intersection of
$\fh$ with the half space determined by $F$ and opposite to $c$.  It is
composed of several $(r-1)$--dimensional domains whose boundaries
compose the boundary $\partial \fh_F$, entirely contained in $F$.
Since we are in a $CAT(0)$ space and $F$ is totally geodesic, the
projection $\pi_F$ is $1$-Lipschitz.  Therefore the volume needed to
fill $\partial \fh_F$ in $F$ is at most $\vol (\fh_F)$.  We denote the new
domain thus obtained by $\widehat{\fh}_F$.  The space has a
bi-combing, therefore $\fh_F \cup \widehat{\fh}_F$ can be filled with
volume at most  $2 \vol (\fh_F) R_F$, where $R_F$ is the minimal radius
$R$ such that the $R$--tubular neighborhood of $F$ contains $\fh_F$.
By hypothesis $\fh$ is $\eta$--round, hence $R_F = O(x)$.

Assume that some of the separating flats $F$ crossed by $\fh$ are at
distance at least $\delta x$ from $c$.  For every such flat we replace
$\fh_F$ by $\widehat{\fh}_F$, and we fill all the $\fh_F \cup
\widehat{\fh}_F$ thus obtained with a volume $O(x^{r+1})$.  Thus, up
to an additional filling volume of order $x^{r+1}$, we may assume that
all the flats crossed by $\fh $ intersect $B(c, \delta x)$.

Without loss of generality we may also assume that the vertex $v_c$ is
in $T_\fh$.  Indeed, if not then let $w$ be the nearest vertex to
$v_c$ in $T_\fh$.  The projection $c'$ of $c$ onto the component
$N(w)$ is the same thing as the projection onto the $r$-separating
flat $F_w$ contained in $N(w)$ and separating it from $c$.  By the
previous, we have that $c'$ is at distance at most $\delta 
x$ and that
consequently $\fh$ is at distance at least $(1-\delta ) x$ from $c'$.
Any filling of $\fh$ that is outside $B\left(c' , \delta x \right)$ is
also outside $B\left(c ,\delta x \right)$: the half space
determined by the flat $F_w$ and opposite to $c$ is convex, therefore
the projection onto it diminishes the distances.  We may thus replace
$c$ by $c'$, at the cost of replacing the parameter $\delta$ by
$\frac{\delta }{1-\delta }$.  Clearly a bounded perturbation of $c'$
allows us to assume that it is contained in only one component, and
that it is a vertex of the triangulation.

Let $F_1,\dots , F_k$ be the separating flats in $N(v_c)$ that are
crossed by $\fh$, and let $e_1,\dots , e_k$ be the corresponding
edges.  Each of the connected components of $T_\fh \setminus \{ v_c
\}$ is a rooted tree $T_i$ with root $v_c$, the first level composed
of only one edge $e_i$ and with $v_c$ removed, for $i\in \{1,2,\dots
,k \}$.  Fix an arbitrary $i\in \{1,2,\dots ,k \}$.  For simplicity we
denote $\fh_{F_i}$ by $\fh_i$.  Let $\ell_i$ be such that the
$(r-1)$-volume of $\fh_i$ is $A\ell_i^{r-1}$.

 Assume first that $\fh_i$ is an $(r-1)$-domain with boundary
 $\partial \fh_i$ an $(r-2)$-hypersurface contained in $F_i$.  Its
 projection onto $F_i$ is a domain of volume at most $A\ell_i^{r-1}$ filling
 $\partial \fh_i$.  Let $c_i$ be the projection of $c$ onto $F_i$.
 According to the previous hypothesis $\dist (c,c_i) \leq \delta
 x$.  In order to avoid $B\left( c, \delta
 x \right)$ when filling a hypersurface in the half
 opposite to $c$ of boundary $F_i$, it suffices to avoid $B\left( c_i,
 \delta x \right)$.  As $\partial \fh_i$ is an
 $(r-2)$-hypersurface in an $r$--dimensional flat, its divergence with
 respect to any point is of the same order as its filling.  Hence,
 there exists an $(r-1)$-domain $\fh_i'$ contained in $F_i$, with
 boundary $\partial \fh_i$, with volume $\leq B\ell_i^{r-1}$, and which
 avoids $B\left( c_i, \delta x \right)$.

 We prove that the hypersurface $\overline{\fh}_i = \fh_i \cup
 \fh_i'$, of volume $\leq C\ell_i^{r-1}$, can be filled outside
 $B\left( c_i, \delta x \right)$ with a domain of volume $\leq D
 \ell_i^{r+n+1}$.  We prove this by induction on the number of
 vertices in the tree $T_i$ that are of valency at least $3$.  If this
 number is $0$ then $\fh_i$ crosses a sequence of consecutive
 separating flats $\widetilde{F}_1,\dots, \widetilde{F}_m$.  Each
 $\widetilde{\fh}_j= \fh_{\widetilde{F}_j}$ has volume $A
 \widetilde{\ell}_j^{r-1}$, for some $\widetilde{\ell}_j$.  We agree
 to denote $\fh_i $ by $\widetilde{\fh}_0$ and correspondingly
 $\ell_i$ by $\widetilde{\ell}_0$.  For each $\widetilde{\fh}_j$ we
 can repeat the argument above and complete it with an
 $\widetilde{\fh}_j'$ contained in $\widetilde{F}_j$ and avoiding
 $B\left( \pi_{\widetilde{F}_j} (c), \delta x \right)$, so that
 $\widetilde{\fh}_j \cup \widetilde{\fh}_j'$ is a hypersurface (or a
 union of hypersurfaces) of volume $\leq C\widetilde{\ell}_i^{r-1}$.
 For $j=m$ the inductive hypothesis on the geometric components
 implies that $\widetilde{\fh}_m \cup \widetilde{\fh}_m'$ can be
 filled outside the ball of center $\pi_{\widetilde{F}_m} (c)$ and of
 radius $\delta x $ by a domain of volume $\leq D
 \widetilde{\ell}_m^{r+n}$.  Indeed, if $\widetilde{\ell}_i \leq
 \varepsilon x$ for $\varepsilon$ small enough then,  
 by the relation between the filling radius and filling function 
 established in \cite{BD:HigherFilling1}, it follows that 
 the ordinary filling of
 $\widetilde{\fh}_m \cup \widetilde{\fh}_m'$ already avoids a ball of
 radius $\delta x $ centered in $\pi_{\widetilde{F}_m} (c)$, while if
 $\widetilde{\ell}_i \geq \varepsilon x$ then one can apply the usual
 estimate of the divergence in $\widetilde{N}_n$ and obtain that
 $\widetilde{\fh}_m \cup \widetilde{\fh}_m'$ can be filled outside
 $B\left( \pi_{\widetilde{F}_m} (c) , \delta x \right)$ with a volume
 $\preceq x^{r+n} \preceq \widetilde{\ell}_m^{r+n}$.  For $j<m$ one
 considers the hypersurface composed of $\widetilde{\fh}_j'$ and
 $\widetilde{\fh}_j$, with $\widetilde{\fh}_{j+1}$ replaced by
 $\widetilde{\fh}_{j+1}'$.  Again the inductive hypothesis implies
 that this domain can be filled outside $B\left( c, \delta x \right)$
 with a volume $\leq D \widetilde{\ell}_j^{r+n}$.  We thus obtain a
 filling of $\overline{\fh}_i$ of volume $\leq D \sum_{j=1}^m
 \widetilde{\ell}_j^{r+n}$.

 \begin{figure}[h]
 \includegraphics[width=0.7\textwidth]{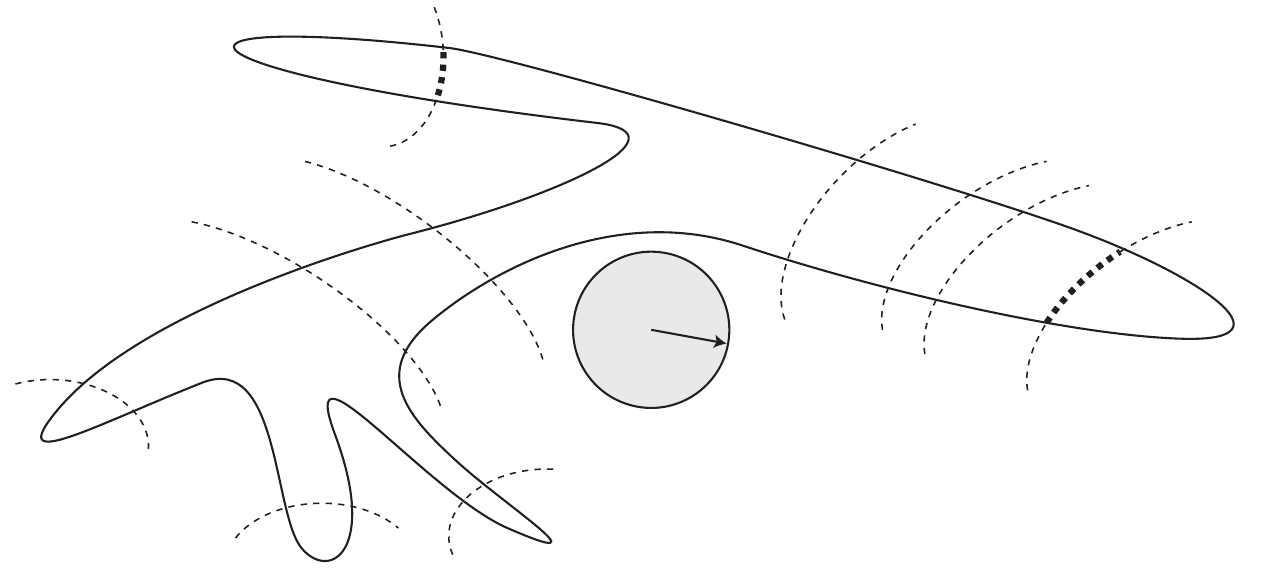}
 \put  (0, 45) {$\widetilde{\fh}_{m}$}
 \put   (-10, 66){$\widetilde{F}_{m}$}
 \put   (-30, 50){$\widetilde{\fh}'_{m}$}
 \put   (-55, 58){$\ldots$}
 \put   (-32, 71){$\widetilde{F}_{2}$}
 \put   (-42, 83){$\widetilde{F}_{1}$}
 \put   (-65, 88){$F_{2}$}
 \put   (-118, 50){$\delta_{x}$}
 \put   (-128, 50){$c$}
 \put  (-160, 95)  {$\widetilde{\fh}'_{1}$}
 \put  (-215, 95)  {$\widetilde{\fh}_{1}$}
 \put  (-170, 72){$\widetilde{F}_{3}$}
 \put   (-80, 60){$\ldots$}
  \put  (-162, 35){$F'$}
 \put  (-140, 15){$F''_{3}$}
 \put  (-220, 5){$F''_{2}$}
 \put  (-263, 35){$F''_{1}$}
 \caption{Set-up for the proof of the divergence upper bound 
 in $\widetilde{N}_{n}$}\label{fig:UpperboundDiv}
 \end{figure}
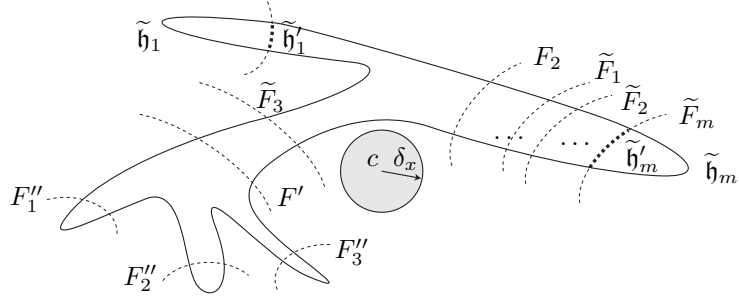

  Note that $\widetilde{\ell}_1,\dots, \widetilde{\ell}_m$ is a
  decreasing sequence of numbers that are at least
  $\frac{1}{A^{1/(r-1)}}$.  If for some $j$ we have
  $\widetilde{\ell}_j=\widetilde{\ell}_{j+1}$ then the part of $\fh$ in
  between $\widetilde{F}_j$ and $\widetilde{F}_{j+1}$ does not contain
  any
  $k$-simplex.  Since this separates $\fh$, by filling this part with a
  domain of volume zero we can split $\fh$ into two hypersurfaces
  $\fh_1,\fh_2$ whose volume adds to that of $\fh$, and argue for
  $\fh_1,\fh_2$ instead of $\fh$.

  Thus without loss of generality we may assume that $A\widetilde{\ell}_j^{r-1}\geq A\widetilde{\ell}_{j+1}^{r-1} +1$ for every $j$, whence $\widetilde{\ell}_j\geq \widetilde{\ell}_{j+1} + \frac{1}{A (r-1)}$. We deduce that
  $$
  \sum_{j=1}^m \widetilde{\ell}_j^{r+n} \leq 2A (r-1) \int_0^{\widetilde{\ell}_0} x^{r+n}\, \mathrm{d}x\, ,
  $$
   and conclude that $\overline{\fh}_i$ can be filled with a volume $\preceq \ell_i^{r+n+1}$.

  Assume that the required estimate holds when $T_i$ has at most
  $a$ vertices of valence $\geq 3$ and now assume that $T_i$ has $a+1$
  vertices with valence $\geq 3$.  Let $v$ be such a vertex that is
  farthest from $v_c$, let $e'$ be the edge adjacent to it that is
  nearest $v_c$, and $e_1,\dots , e_r$ the other edges.  The
  corresponding separating flats are $F', F_1,....F_r$, their intersections with $\fh$ are the
  hypersurfaces $\fh' , \fh_1,\dots , \fh_r$, and the respective volumes of these hypersurfaces are $A(\ell')^{r-1},A(\ell_1)^{r-1}, \ldots, A(\ell_r)^{r-1}$.  With an argument as above, we
  reduce to the case when $\fh_i \subseteq F_i$, at the cost of some
  volume $D\sum_{i=1}^r \ell_i^{r+n+1}$.  The hypersurface thus
  obtained can be filled in $N(v)$ outside the required ball with a
  volume $\leq D(\ell' )^{r+n} $.  On the whole we obtain that we can
  assume that $\fh' \subset F'$ at the cost of some volume $\leq
  D(\ell' )^{r+n+1}$.  With this change, we now can apply the inductive
  argument.

  If $\fh_i$ is composed of several $(r-1)$-domains with boundaries in $F_i$ then a similar argument, repeated for each of these domains, implies that one can assume that $\fh_i \subseteq F_i$, at the cost of a volume $\leq D(\ell_i )^{r+n+1}$.

  The inductive hypothesis applied again in $N(v_c)$ allows us to conclude
  that to fill $\fh $ outside the required ball the necessary volume is
  $B\ell^{r+n} + \sum_{i=1}^k D(\ell_i )^{r+n+1}$, and thus, in
  particular, this volume is $\preceq
  \ell^{r+n+1}$.
 \endproof


\newcommand{\etalchar}[1]{$^{#1}$}
\def\cprime{$'$} \def\cprime{$'$} \def\cprime{$'$} \def\cprime{$'$}
  \def\cprime{$'$} \def\cprime{$'$}
\providecommand{\bysame}{\leavevmode\hbox to3em{\hrulefill}\thinspace}
\providecommand{\MR}{\relax\ifhmode\unskip\space\fi MR }
\providecommand{\MRhref}[2]{%
  \href{http://www.ams.org/mathscinet-getitem?mr=#1}{#2}
}
\providecommand{\href}[2]{#2}

\end{document}